\documentclass[12pt]{article}
\def\date{28 March 2011}
\usepackage{amsmath, amssymb, amsfonts, theorem}

\newtheorem{proposition}{proposition}[section]
\newtheorem{lemma}[proposition]{Lemma}

\newtheorem{theorem}[proposition]{Theorem}

\newtheorem{conjecture}[proposition]{Conjecture}
\newcommand{\qed}{$\square$\bigskip}

\newcommand{\proof}{{\noindent\bf Proof. }}
{\theorembodyfont{\rmfamily}

}

\def\myclaim#1#2{{\medbreak\noindent\rlap{\rm(#1)}\ignorespaces
 \rightskip20pt
 \hangindent=20pt\hskip20pt{ \ignorespaces\sl#2}\smallskip}}

\begin{document}
\font\smallrm=cmr8




\baselineskip=12pt
\phantom{a}\vskip .25in
\centerline{{\bf  SUB-EXPONENTIALLY MANY $3$-COLORINGS}}
\centerline{{\bf  OF TRIANGLE-FREE PLANAR GRAPHS}}
\vskip.4in
\centerline{{\bf Arash Asadi}%
\footnote{\texttt{aasadi@math.gatech.edu}.}}
\smallskip
\centerline{School of Mathematics}
\centerline{Georgia Institute of Technology}
\centerline{Atlanta, Georgia  30332-0160, USA}
\medskip
\centerline{{\bf Zden\v ek Dvo\v r\'ak}%
\footnote{\texttt{rakdver@kam.mff.cuni.cz}.  Supported by Institute for Theoretical Computer Science (ITI), project 1M0021620808 of Ministry of Education of Czech Republic, and by project GA201/09/0197 (Graph colorings and flows: structure and applications) of Czech Science Foundation.}}
\smallskip
\centerline{Department of Applied Mathematics}
\centerline{Charles University}
\centerline{Malostransk\'e n\'am.~25}
\centerline{118 00 Prague 1, Czech Republic}

\medskip
\centerline{{\bf Luke Postle}%
\footnote{\texttt{ljpostle@math.gatech.edu}. 
Partially supported by an NSF Graduate Research Fellowship.}}
\centerline{{\bf Robin Thomas}%
\footnote{\texttt{thomas@math.gatech.edu}. Partially supported by NSF under
Grant No.~DMS-0701077.}}
\smallskip
\centerline{School of Mathematics}
\centerline{Georgia Institute of Technology}
\centerline{Atlanta, Georgia  30332-0160, USA}

\vskip 1in \centerline{\bf ABSTRACT}
\bigskip

{
\parshape=1.0truein 5.5truein
\noindent
Thomassen conjectured that every triangle-free planar graph on $n$ vertices
has exponentially many {$3$-colorings}, and proved that it has at least
$2^{n^{1/12}/20000}$ distinct $3$-colorings. 
We show that it has at least
$2^{\sqrt{n/212}}$ distinct $3$-colorings.
}

\vfill \baselineskip 11pt \noindent April 2010, revised \date.
An extended abstract of this paper appeared in~\cite{AsaPosThosubabstr}.
\vfil\eject
\baselineskip 18pt

\section{Introduction}
All graphs in this paper are finite, and have no loops or multiple
edges. Our terminology is standard, and may be found
in~\cite{Bollobas} or~\cite{Diestel}.
 In particular, cycles and paths have no repeated vertices.
The following is a well-known theorem of Gr\"otzsch~\cite{Grotzsch}. 

\begin{theorem}
\label{thm:grotzsch}
Every triangle-free planar graph is $3$-colorable.
\end{theorem}

Theorem~\ref{thm:grotzsch} has been the subject of extensive research.
Thomassen~\cite{Thomassen_projective} gave several short 
proofs~\cite{Thomassen_projective, Thomassen_3listgirth5,
  Thomassen_short}
of Gr\"otzsch's theorem and extended it to 
projective planar and toroidal graphs. 
The theorem does not extend verbatim to any non-planar surface,
but Thomassen proved that every graph of girth at least five
embedded in the projective plane or the torus is $3$-colorable.
Gimbel and Thomassen~\cite{GimTho} found an elegant characterization 
of $3$-colorability for triangle-free projective planar graphs.
There does not seem to be a corresponding counterpart for other surfaces,
but Kr\'al' and Thomas~\cite{KraTho} found a characterization of 
$3$-colorability for toroidal and Klein bottle graphs that are embedded
with all faces even.
It was an open question for a while whether a $3$-coloring of a triangle-free
planar graph can be found in linear time. 
First Kowalik~\cite{Kowalik}
 designed an almost linear time algorithm, and then a linear-time
algorithm was found by Dvo\v{r}\'ak, Kawarabayashi and 
Thomas in~\cite{ZKTlinear}.
For a general surface $\Sigma$, Dvo\v{r}\'ak, Kr\'al' and 
Thomas~\cite{DvoKraTho} found a linear-time algorithm to decide whether
a triangle-free graph in $\Sigma$ is $3$-colorable.

In this paper we study how many $3$-colorings 
a triangle-free planar graph must have.
Thomassen conjectured in~\cite{Thomassen_manycoloring} that exponentially many:

\begin{conjecture}
\label{expconj} 
There exists an absolute constant $c>0$, such that
if $G$ is a triangle-free planar graph on $n$ vertices, then $G$
has at least $2^{cn}$ distinct $3$-colorings.
\end{conjecture}

\noindent
Thomassen gave a short proof of this conjecture under the additional
hypothesis that $G$ has girth at least five.
We use that argument in Lemma~\ref{exp_antichain} below;
Thomassen's original proof may be recovered by taking $\cal F$ to be the
set of all facial cycles.
Thomassen~\cite{Thomassen_manycoloring} then extended this
result by showing that every planar graph of girth at least five
has exponentially many list-colorings for every list assignment that
gives each vertex a list of size at least three.
For triangle-free graphs Thomassen~\cite{Thomassen_manycoloring} 
proved a weaker version of Conjecture~\ref{expconj}, namely that every 
triangle-free planar graph on $n$ vertices has at least $2^{n^{1/12}/20000}$
distinct $3$-colorings.
Our main result is the following improvement.

\begin{theorem}
\label{maintheorem}
Every triangle-free planar graph on $n$ vertices has
at least $2^{\sqrt{n/212}}$ distinct $3$-colorings.
\end{theorem}

In closely related work Thomassen~\cite{Thomassen_choosableexp} 
proved that every (not necessarily triangle-free) planar graph has 
exponentially many list colorings provided every vertex has at
least five available colors.

Our paper is organized as follows. In the next section we investigate
non-crossing families of $5$-cycles, and reduce Theorem~\ref{maintheorem}
to Lemma~\ref{exp_chain}, which states that if a triangle-free
planar graph has $k$ nested $5$-cycles, then it has at least
$2^{k/7}$ $3$-colorings.
The rest of the paper is devoted to a proof of Lemma~\ref{exp_chain},
which we complete in Section~\ref{sec:chains}.
In Section~\ref{sec:matrix} we prove an auxiliary result stating
that some entries in the product of certain matrices grow exponentially
in the number of matrices.

We end this section by stating a useful theorem of
Thomassen~\cite{Thomassen_projective}.

\begin{theorem} \label{ThomExtension45}
Let $G$ be a triangle-free plane graph with facial cycle $C$ of length at most five. Then every $3$-coloring of $C$ extends to a $3$-coloring of $G$.
\end{theorem}

We would like to acknowledge that an extended abstract 
of this paper appeared in~\cite{AsaPosThosubabstr}.

\section{Laminar Families of $5$-Cycles}
First we define some terminology. Let $A$ and $B$ be two subsets
of $\mathbb{R}^2$. We say that $A$ and $B$ \emph{cross} if $A\cap
B$, $A\cap B^c$, $A^c \cap B$, $A^c\cap B^c$ are all non-null.
Then we say that a family $\mathcal{F}$ of subsets of
$\mathbb{R}^2$ is \emph{laminar} if for every two sets $A,B\in
\mathcal{F}$, $A$ and $B$ do not cross. Now let $G$ be a plane
graph and $C$ be a cycle in $G$. Then we let $Int(C)$ denote the
bounded component of $\mathbb{R}^2-C$ and $Ext(C)$ denote the
unbounded component of $\mathbb{R}^2-C$. Now we say that a family $\mathcal{F}$
of cycles of $G$ is \emph{laminar} if the
corresponding family of sets $\{Int(C): C\in \mathcal{F}\}$
is laminar. 
We call a family $\mathcal{F}$ of cycles an
\emph{antichain} if $Int(C_1)\cap Int(C_2) = \emptyset$
for every distinct $C_1, C_2\in\mathcal{F}$, and we call it
a \emph{chain} if for every two cycles $C_1,
C_2\in\mathcal{F}$, either $Int(C_1)\subseteq Int(C_2)$ or
$Int(C_2)\subseteq Int(C_1)$.

Let $G$ be a triangle-free plane graph, and let $v\in V(G)$. We define $G_v$ to
be the graph obtained from $G$ by deleting $v$, identifying all
the neighbors of $v$ to one vertex, and deleting resulting
parallel edges. We also let $D_k(G)$ denote
 the set of vertices of $G$ with degree at most $k$.

\begin{lemma}{ \label{bd_or_five} If $G$ is a triangle-free plane graph and $k \geq 0$ is an integer, then either
\begin{itemize}
\item[(i)] there exists $v\in D_k(G)$ such that $G_v$ is
triangle-free
 or,
\item[(ii)] there exists a laminar family $\mathcal{F}$ of $5$-cycles
such that every $v\in D_k(G)$ belongs to some member of $\mathcal{F}$.
\end{itemize}
}\end{lemma}

\proof
 We proceed by induction on the number of vertices of $G$. Suppose
condition (i) does not hold. Notice that if $v\in V(G)$ and $G_v$
is not triangle-free, this implies, since $G$ is triangle-free,
that $v$ is in a $5$-cycle in $G$. Hence if condition (i) does not
hold, every $v\in D_k(G)$ must be in a $5$-cycle in $G$.

Now suppose there does not exist a separating $5$-cycle in $G$.
Then we let $\mathcal{F}$ be the set of all $5$-cycles in $G$. The
second condition then holds since the absence of separating cycles
implies that $\mathcal{F}$ is laminar.

Thus we may assume that there exists a $5$-cycle $C$ that
separates $G$ into two triangle-free plane graphs $G_1$ and
$G_2$, where both $G_1$ and $G_2$ include $C$. 
By induction, the lemma holds for $G_1$ and $G_2$. Suppose
that both $G_1$ and $G_2$ satisfy condition (ii) with laminar
families $\mathcal{F}_1$ and $\mathcal{F}_2$, respectively. Then
let $\mathcal{F}=\mathcal{F}_1 \cup \mathcal{F}_2$. Note that 
$\mathcal{F}$ is laminar. 
Now $G$ satisfies condition (ii) since every
$v\in D_k(G)$ is contained in either $D_k(G_1)$ or $D_k(G_2)$.
Thus we may assume without loss of generality that $G_1$ satisfies
condition (i). That is, there exists $v\in D_k(G_1)$ such that
$(G_1)_v$ is triangle-free. This implies that $v$ is not in a
$5$-cycle in $G_1$. 
In particular, $v\not\in V(C)$, and hence $v\in D_k(G)$.
Yet since $G_v$ is not triangle-free
by assumption, $v$ must be in a $5$-cycle in $G$, say $C'$. It follows that $C'$
intersects $C$.  Let $P$ be the maximal subpath of $C'$ containing $v$ such that no
internal vertex of $P$ belongs to $C$.  Since $G$ is triangle-free, $P$ has length
$t\in\{2,3\}$ and the endvertices of $P$ are joined by a path of length $5-t$ contained in $C$.
Hence $v$ is in a 5-cycle in $G_1$, a contradiction.~\qed

\begin{lemma}{ \label{bd_anti_chain} 
If $G$ is a triangle-free plane graph on $n$ vertices,
then either
\begin{itemize}
\item[(i)] there exists $v\in D_k(G)$ such that $G_v$ is
triangle-free, or
\item[(ii)] $G$ has an antichain $\mathcal{F}$ of $5$-cycles such
that $|\mathcal{F}| \ge \sqrt{{6(k-3)n\over35(k-1)}}$, or
\item[(iii)] $G$ has a chain $\mathcal{F}$ of $5$-cycles such that
$|\mathcal{F}| \ge \sqrt{{7(k-3)n\over30(k-1)}}$.
\end{itemize}
}\end{lemma}

\proof
Since $G$ is  triangle-free and planar, it satisfies $2|V(G)| \geq |E(G)|$. 
We may assume that (i) does not hold and hence every vertex of $G$ has degree at least two. It follows that
$$4|V(G)| \geq 2 |E(G)|=\sum_{v \in V(G)}deg(v) \geq (k+1)\left(|V(G)|-|D_k(G)|\right)+2|D_k(G)|,$$
and hence $|D_k(G)| \geq \frac{k-3}{k-1}|V(G)|$. 
Since (i) does not hold, we deduce
from Lemma~\ref{bd_or_five} that
there exists a laminar family of $5$-cycles $\mathcal{G}$ of size
at least $|D_k(G)|/5\ge\frac{k-3}{5(k-1)}n$. 
By Dilworth's theorem applied to the partial order on $\mathcal{G}$
defined by $Int(C_1) \subseteq Int(C_2)$ we deduce
that $\mathcal{G}$ has either an antichain of size at least
$\sqrt{6|\mathcal{G}|/7}$, in which case condition (ii) holds, or
a chain of size at least $\sqrt{7|\mathcal{G}|/6}$, in
which case condition (iii) holds.~\qed

\begin{lemma} 
\label{exp_antichain} 
Let $G$ be a triangle-free plane graph. If $G$ has an
antichain $\mathcal{F}$ of $5$-cycles, then $G$ has at least
$2^{|\mathcal{F}|/6}$ distinct 3-colorings.
\end{lemma}

\proof
Let $G'$ be obtained from $G$ by deleting the vertices in
$\bigcup_{C\in\mathcal{F}}Int(C)$. Now $G'$ has at least
$|\mathcal{F}|$ facial $5$-cycles. 
By Euler's formula $|E(G')|\le2|V(G')|-|\mathcal{F}|/2$.
By Theorem~\ref{thm:grotzsch} the graph $G'$ has a 3-coloring $\Phi$. 
For  $i,j\in\{1,2,3\}$ with $i<j$  let $G_{ij}$ denote the subgraph of $G$
induced by the vertices colored $i$ or $j$. Since 
$\sum_{i<j}(|V(G_{ij}|-|E(G_{ij}|)=2|V(G')|-|E(G')|\ge |\mathcal{F}|/2$, 
there exist $i,j\in\{1,2,3\}$ such that $i<j$ and $G_{ij}$ 
has at least $|\mathcal{F}|/6$ components. But then
there are at least $2^{|\mathcal{F}|/6}$ distinct 3-colorings of $G'$ since
switching the colors on any subset of the components of $G_{ij}$
gives rise to a distinct coloring of $G'$.
Furthermore, every 3-coloring of $G'$
extends to a 3-coloring of $G$ by Theorem~\ref{ThomExtension45}.~\qed

\begin{lemma}
\label{exp_chain}
Let $G$ be a triangle-free plane graph. If $G$ has a chain
$\mathcal{F}$ of $5$-cycles, then $G$ has at least
$2^{|\mathcal{F}|/7}$ distinct $3$-colorings.
\end{lemma}

We will prove Lemma~\ref{exp_chain} in Section~\ref{sec:chains},
but now we deduce the main theorem from it.
\medskip

\noindent
{\bf Proof of Theorem~\ref{maintheorem}, assuming Lemma~\ref{exp_chain}.}
We proceed by induction on the number of vertices. 
If $n\le212$, then the conclusion clearly holds. 
We may therefore assume that $n\ge213$ and that the theorem holds for
all graphs on fewer than $n$ vertices.
If there exists
$v\in D_{213}(G)$ such that the graph $G_v$ 
(defined prior to Lemma~\ref{bd_or_five}) is triangle-free, then by
induction $G_v$ has at least $2^{\sqrt{(n-\deg(v))/212}}$ distinct
$3$-colorings. Hence $G$ has at least $2 \cdot
2^{\sqrt{(n-\deg(v))/212}}$ distinct 3-colorings, which is greater
than $2^{\sqrt{n/212}}$ since $\deg(v) \leq 213$. So we may assume
by Lemma~\ref{bd_anti_chain} applied to $k=213$ 
that $G$ has either an antichain of
$5$-cycles of size at least $\sqrt{36n/212}$, 
in which case the theorem holds by Lemma~\ref{exp_antichain};
 or a chain of $5$-cycles of size at least $\sqrt{49n/212}$,
in which case the theorem holds by Lemma~\ref{exp_chain}.~\qed 


\section{A matrix lemma}
\label{sec:matrix}

Let the matrix $A_0$ be defined by
$$A_0=\left[{\begin{matrix} 1 & 1 & 0 & 0 & 0 \cr 1 & 1 & 0 & 0 & 0 \cr 0 & 0 & 1 & 0 & 0 \cr 0 & 0 & 0 & 1 & 0 \cr 0 & 0 & 0 & 0 & 1 \end{matrix}}\right].$$
Let $A$ and $B$ be two $5\times5$ matrices with non-negative entries.
We say that $A$ \emph{majorizes}  $B$ if every entry in $A$ is greater than or equal to the corresponding entry
of $B$. We say that $A$ \emph{dominates} $B$ if there exist permutation matrices $P,Q$ such that 
$A$ majorizes $PBQ$. We say that $A$ is \emph{dominant} if $A$ dominates $A_0$.
We say that $A$ is \emph{doubling} if, in every row and column of $A$, there are at least two entries with value at least one.\

We denote the vector of all ones by $\bf1$.

\begin{lemma}
\label{MMM}
Let $n\ge1$ be an integer.
If for every $i=1,2,\ldots,n$, 
$M_i$ is either a dominant or a doubling $5\times5$ matrix, then ${\bf1}^TM_{1}M_{2} \cdots M_{n}{\bf1} \ge (3/2)^{n/4}$. 
\end{lemma}

\proof
Let $x = (x_1,x_2,x_3,x_4,x_5)$ be a vector of non-negative integers. Let $s_k(x)$ denote the sum of the $k$ smallest entries of $x$, and let $S(x) = s_1(x)s_2(x)s_4(x)s_5(x)$. If $M$ majorizes $N$, then $s_k(Mx) \ge s_k(Nx)$ for every $k\in \{1,2,3,4,5\}$, and hence $S(Mx)\ge S(Nx)$. If $P$ is a permutation matrix, then $s_k(Px)=s_k(x)$ for every $k\in \{1,2,3,4,5\}$, and hence $S(Px)=S(x)$. Consequently, if $M$ dominates $N$, then $s_k(Mx)\ge s_k(Nx)$ for every $k\in\{1,2,3,4,5\}$, and hence $S(Mx)\ge S(Nx)$. 

We claim that if $M$ is dominant, then $S(Mx) \ge 3 S(x)/2$
for every vector $x$ of non-negative integers. 
>From above, it is sufficient to prove that $S(A_0x)\ge 3S(x)/2$
for every $x$. 
As $A_0$ dominates the identity matrix, $s_k(A_0x)\ge s_k(x)$ for every $k\in\{1,2,3,4,5\}$. Without loss of generality suppose that $x_1\le x_2$ and $x_3 \le x_4 \le x_5$.  
\begin{itemize}
\item If $x_1 + x_2 \le x_3$, then $s_1(A_0x) = x_1 + x_2 \geq 2s_1(x)$.
\item If $x_3 < x_1 + x_2 \le x_4$, then $s_2(A_0x) = x_3 + (x_1+x_2)$ 
and $s_2(x)= x_1 + \min\{x_3,x_2\}$. 
Thus $s_2(A_0x) \ge 3 s_2(x)/2$.
\item If $x_4 < x_1 + x_2 \le x_5$, then $s_4(A_0x) = 2x_1+2x_2+x_3+x_4$, $s_4(x) =x_1+x_2+x_3+x_4$, $s_2(A_0x)=x_3+x_4$, and $s_2(x)= \min\{x_1+x_2,x_1+x_3,x_3+x_4\}$. 
If $x_3+x_4\le x_1+x_2$, then $s_4(A_0x)\ge3s_4(x)/2$, as desired.  On the other hand, if $x_3+x_4\ge x_1+x_2$,
then note that $s_4(A_0x)\ge 3s_2(x)$ and $s_2(A_0x)\ge s_4(x)/2$; 
hence $$\frac{s_2(A_0x)s_4(A_0x)}{s_2(x)s_4(x)}\ge {3\over2}\,.$$
\item If $x_5 < x_1 + x_2$, then 
$s_4(A_0x)=x_1+x_2+x_3+x_4+x_5 \ge 5s_4(x)/4$ and 
$s_5(A_0x) = s_5(x)+x_1+x_2 > 5s_5(x)/4$. 
Hence, $s_4(A_0x)s_5(A_0x) \ge 25s_4(x)s_5(x)/16$ and 
$S(A_0x) \ge 25S(x)/16\ge3S(x)/2$.
\end{itemize}
This proves the claim.

We claim that if $M$ is doubling, then $S(Mx) \ge 10 S(x)$. 
Recall that, by the definition of doubling matrix, each row and each column of 
$M$ contains at least two entries greater or equal to one. Thus each entry of $Mx$ is a sum containing at least 
two entries of $x$, and each entry of $x$ appears in at least two entries 
of $Mx$. It follows that $s_1(Mx) \ge 2s_1(x)$, $s_2(Mx) \ge 2s_2(x)$ 
and $s_5(Mx) \ge 2s_5(x)$.  
Furthermore, a sum of any four entries of $Mx$ contains all entries of $x$, 
and hence $s_4(Mx) \ge x_1+x_2+x_3+x_4+x_5 \ge 5 s_4(x)/4$. 
Thus $S(Mx) \ge 10 S(x)$ and the claim is proved.

Let $x_n={\bf1}^TM_{1}M_{2} \cdots M_{n}$.
Note that $s_k(x_0)=k$, thus $S(x_0)=40$.
As $M_i$ is dominant or doubling for every $i$, $S(x_n)\ge 40 (3/2)^n$.
Finally, $x_n{\bf1} = s_5(x_n) \ge S(x_n)^{1/4}$ and the lemma follows.
~\qed

\section{ Chains of $5$-Cycles }
\label{sec:chains}
In order to prove Lemma~2.4, we will first characterize how the
3-colorings of an outer $5$-cycle of a plane graph $G$ extend to 
the 3-colorings of another $5$-cycle. 
If $C$ is a $5$-cycle in a graph $G$ and $\Phi$ a
3-coloring of $C$, then there exists a unique vertex $v\in V(C)$
such that $v$ is the only vertex of $C$ colored $\Phi(v)$.
We call such a vertex
the \emph{special vertex of $C$ for $\Phi$}. 
Let $e$ be the edge of $C$ opposite the special vertex of $C$ for $\Phi$.
We call such an edge the \emph{special edge of $C$ for $\Phi$}.

Let $G$ be a triangle-free plane graph and $C_1, C_2$ be $5$-cycles in $G$
such that $C_1\ne C_2$ and $Int(C_2)\subseteq Int(C_1)$. 
Let $C_1:= u_1\ldots u_5, C_2:=v_1\ldots v_5$. Then we
define a \emph{color transition matrix} $M$ of $G$ with respect to
$C_1$ and $C_2$ as follows. 
Let $G'$ be the subgraph of $G$ consisting of all
the vertices and edges of $G$ drawn in the closed
annulus bounded by $C_1\cup C_2$.
We let $M_{ij}$ equal one sixth the number of 3-colorings
$\Phi$ of $G'$ such that $u_i$ is the special vertex of $C_1$ for $\Phi$
and $v_j$ is the special vertex of $C_2$ for $\Phi$. Note that $A_0$ is a color transition matrix
of a graph $G$ when $G=C_1\cup C_2$ and $C_1$ and $C_2$ have four vertices in common.

The following lemma is straightforward.

\begin{lemma}{\label{multiplication}
Let $n\ge2$ be an integer,
let $G$ be a triangle-free graph and $\mathcal{F}=\{C_1, \dots,
C_n\}$ be a family of $5$-cycles such that  ${Int(C_i) \supseteq
Int(C_j)}$ if  $1 \leq i < j \leq n$. Let $M_i$ be a
color transition matrix of $G$ with respect to $C_i$ and
$C_{i+1}$. Then $M_1M_{2}\ldots M_{n-1}$ is a color transition matrix of $G$ with respect to
$C_1$ and $C_n$.}
\end{lemma}



To prove our next lemma, we need the following theorem which follows by 
combining Theorems 1 and 2 of Aksionov~\cite{Aksionov}:

\begin{theorem}
\label{Aksionov}
Let $G$ be a plane graph with facial $5$-cycle $C$ and exactly one triangle $T$ which is facial. Let $\Phi$ be a
$3$-coloring of $C$ and $e$ be the special edge of $C$ for $\Phi$. Then $\Phi$ does not extend to a $3$-coloring
of $G$ if and only if $e$ is an edge of $T$ and $G$ has a subgraph $H$,
where every face of $H$ has length four except for the faces bounded by
$C$ and $T$.
\end{theorem}

\begin{lemma}
\label{transition} 
{
Let $G$ be a triangle-free plane graph and $C_1, C_2$ be two distinct
$5$-cycles in $G$. Every color transition matrix of $G$ with respect to $C_1$ and
$C_2$ is either dominant or doubling.
}
\end{lemma}

\proof
Let us assume for a contradiction that the lemma is false, and choose
a counterexample $G$ with cycles $C_1$ and $C_2$ with $|V(G)|$ minimum.
Let $M$ be a color transition matrix of $G$ with respect to $C_1$ and $C_2$
that is neither dominant nor doubling.
Let $C_1:=u_1u_2u_3u_4u_5$ and $C_2:=v_1v_2v_3v_4v_5$. 

\myclaim{1}{Every $4$-cycle $C$ in $G$ separates $C_1$ from $C_2$.}

\noindent
To prove (1), suppose for a contradiction that a $4$-cycle $C=x_1x_2x_3x_4$ does not separate $C_1$ from $C_2$.
First suppose that $C$ is not facial. Then some component $J$ of $G\backslash V(C)$ is disjoint from
$C_1\cup C_2$, and hence every $3$-coloring of $G\backslash V(J)$ extends to
$G$ by Theorem~\ref{ThomExtension45}. Thus $M$ dominates every color transition matrix of $G \backslash V(J)$ with respect
to $C_1$ and $C_2$. Hence, by the minimality of $G$, $M$ is either dominant or doubling,
 a contradiction. 

So we may assume that $C$ is facial. Let $G_1$ be the graph
obtained from $G$ identifying $x_1$ and $x_3$ and let $G_2$ be the
graph obtained from $G$ by identifying $x_2$ and $x_4$. At least
one of the graphs $G_1$, $G_2$ is a triangle-free
plane graph. From the symmetry we may assume that $G_1$ is triangle-free. 
Let $C_1'$ and $C_2'$ be the $5$-cycles in $G_1$ corresponding to $C_1$ and $C_2$.
As every $3$-coloring of $G_1$ extends to a $3$-coloring of $G$,
$M$ dominates the color transition matrix of $G_1$ with respect to $C_1',C_2'$.
If $C_1'$ is distinct from $C_2'$, then $M$ is dominant or doubling by the minimality of $G$, a contradiction.
So we may assume that $C_1'=C_2'$. Thus $G=C_1\cup C_2$ where $C_1$ and $C_2$ have four vertices in common. As noted earlier, every color transition matrix of $G$ with respect to $C_1,C_2$ dominates $A_0$. 
Hence $M$ dominates $A_0$ and $M$ is dominant, a contradiction.
This proves (1).
\medskip

\myclaim{2}{Let $\Phi$ be a $3$-coloring of $C_1$. If $M$ is not dominant,
then there exist two $3$-colorings of $G$, $\Phi_1$ and $\Phi_2$, extending $\Phi$ such that
the special vertex of $C_2$ for $\Phi_1$ is distinct from the special vertex of $C_2$ for $\Phi_2$.}

\noindent
To prove (2) we note that $\Phi$ extends to a coloring $\Phi_1$ of $G$
by Lemma~\ref{ThomExtension45}.
We may assume without loss of generality that $e_1=u_1u_2$ is the special edge of $C_1$ for $\Phi$. 
We may also assume without loss of generality that $e_2=v_1v_2$ is the special edge of $C_2$ for $\Phi_1$ and hence that $v_4$ is the special vertex of $C_2$ for $\Phi$. Let $G_1$ be obtained from $G$ by adding the edge $v_1v_3$ and $G_2$ be obtained from $G$ by adding
the edge $v_2v_5$.

We claim that if $v_1v_3$ is not a chord of $C_1$ and $G_1$ contains exactly one triangle $T_1=v_1v_2v_3$, then (2) follows.
If $\Phi$ extends to a $3$-coloring $\Phi_2$ of $G_1$, then $\Phi_2$ is also a $3$-coloring of $G$ where the special vertex of $C_2$ for $\Phi_2$ is one of $v_1,v_2,v_5$ and hence distinct from $v_4$, as desired. So we may assume by Theorem~\ref{Aksionov} that $G$ contains a subgraph $H$ where the faces of $H$ have length four except for $C_1$ and $T_1$. Since $v_1v_3$ is not a chord of $C_1$, $H$ has a face $f$ of length four that is not incident with $v_1v_3$. But then $f$ is bounded by a $4$-cycle in $G$ that does not separate $C_1$ from $C_2$, contradicting (1).
By symmetry, if $v_2v_5$ is not a chord of $C_1$ and $G_2$ contains exactly one triangle $T_2=v_1v_2v_5$, then (2) follows.

So suppose that $G_1$ contains more than one triangle. Hence there is a vertex $z\not\in V(C_2)$ adjacent to $v_1$ and $v_3$. By (1), the $4$-cycle $v_1v_2v_3z$ separates $C_1$ from $C_2$. So $v_4$ and $v_5$ are not in $C_1$.
In this case, $G_2$ contains exactly one triangle $T=v_1v_2v_5$ as $G$ is triangle-free. Moreover, $v_2v_5$ is not a chord of $C_1$ as $v_5$ is not in $C_1$. As noted above, (2) now follows.

So we may suppose that $G_1$ contains exactly one triangle $T_1=v_1v_2v_3$ and that $G_2$ contains exactly one triangle $T_2=v_1v_2v_5$. If $v_1v_3$ is not a chord of $C_1$, (2) follows as above. So $v_1v_3$ is a chord of $C_1$. Similarly, $v_2v_5$ must also be a chord of $C_1$. Hence $v_1,v_2,v_3,v_5\in C_1$. As $C_1$ is distinct from $C_2$, $v_4\not\in C_1$. Let $z$ be the remaining vertex of $C_1$. Now $v_3v_4v_5z$ is a $4$-cycle in $G$ that does not separate $C_1$ from $C_2$, contradicting (1). This proves (2).

\medskip

We will now prove that if $M$ is not dominant, then $M$ is doubling; that is,
there are at least two nonzero entries in every row and column of $M$.
Notice that a row $i$ of $M$ will contain at least two nonzero entries if and only if every coloring of $C_1$
with special vertex $u_i$ extends to at least two colorings of $G$, where the induced colorings of $C_2$ have distinct special vertices.
This property follows from (2). By the symmetry between $C_1$ and $C_2$, 
(2) also implies that every column contains two nonzero entries.~\qed

\noindent
{\bf Proof of Lemma~\ref{exp_chain}.}
Let $n=|\mathcal{F}|$. If $n=1$, the lemma follows from Theorem~\ref{thm:grotzsch}. So we may assume that $n\ge2$. 
Let $C_1, C_2, \ldots, C_{n}$ be the elements of $\mathcal{F}$ such that ${Int(C_i) \supseteq Int(C_j)}$
if and only if $1 \leq i < j \leq n$. For $i=1,2,\ldots,n-1$ 
let $M_i$ be a color transition matrix of $G$ with respect to $C_i,
C_{i+1}$. 
Lemma~\ref{multiplication} implies that $M=M_{1}M_{2} \ldots M_{n-1}$
is a color transition matrix of $G$ with
respect to $C_1,C_n$. 
Hence the number of $3$-colorings of $G$ is at least six times
${\bf 1}^TM{\bf 1}$. 
For all $1 \leq i \leq n-1$, Lemma~\ref{transition} implies that $M_i$ is either dominant or doubling.
It follows from Lemma~\ref{MMM} that the number of $3$-colorings of $G$ 
is at least $6\cdot(3/2)^{(n-1)/4} \ge (3/2)^{n/4}\ge2^{n/7}$, as desired.~\qed

%

\baselineskip 11pt
\vfill
\noindent
This material is based upon work supported by the National Science Foundation.
Any opinions, findings, and conclusions or
recommendations expressed in this material are those of the authors and do
not necessarily reflect the views of the National Science Foundation.
\eject

\end{document}